\documentclass[11pt]{article}
\usepackage{fullpage}
\usepackage{authblk}
\usepackage{natbib}
\usepackage{algorithm}
\usepackage{algpseudocode}
\usepackage{amsmath,amsthm,amsfonts}
\usepackage{amsmath}
\usepackage[subnum]{cases}
\usepackage{tcolorbox}
\usepackage{mathrsfs}
\usepackage{amssymb}
\usepackage{color}
\usepackage{mathrsfs}
\usepackage{empheq}
\usepackage{enumitem}
\usepackage{bm}
\usepackage{multirow}
\usepackage{booktabs}
\usepackage{makecell}
\usepackage{graphicx}
\usepackage{subcaption}
\usepackage{comment}
\usepackage{cases}
\usepackage{appendix}
\usepackage{tikz}
\usetikzlibrary{arrows,shapes}
\usepackage[colorlinks= true, linkcolor=red, citecolor=blue, urlcolor=black]{hyperref}
\usepackage{cleveref}
\usepackage{marginnote}
\usepackage{diagbox} 
\usepackage{pgfgantt}

\numberwithin{equation}{section}

\theoremstyle{definition}

\newtheorem{theorem}{Theorem}[section]
\newtheorem{lemma}[theorem]{Lemma}

\newtheorem{remark}[theorem]{Remark}

\newtheorem{defn}[theorem]{Definition}

\DeclareMathOperator*{\argmax}{argmax}
\DeclareMathOperator*{\argmin}{argmin}

\usepackage{booktabs}
\usepackage{tikz}
\usepackage{pgfplots}
\usetikzlibrary{calc,intersections}

\usetikzlibrary{shapes.geometric, arrows, positioning}

\tikzset{
	mybox/.style  = {draw, rectangle, minimum width=4cm, minimum height=0.8cm, text centered, text width=4.4cm,   
		font=\normalsize},
	box/.style  = {draw, rectangle, minimum width=2.0cm, minimum height=0.6cm, text centered, text width=3.0cm,   
		font=\normalsize},
	myarrow/.style = {line width=0.2pt, draw=black, -triangle 60, postaction={draw, line width=0.2pt, shorten >=10pt,-}}
}

\tikzstyle{arrow} = [->, >=stealth, -triangle 60]

\allowdisplaybreaks

\makeatletter
\newcommand{\leqnomode}{\tagsleft@true}
\newcommand{\reqnomode}{\tagsleft@false}
\makeatother

\begin{document}
%\maketitle

\title{A Lyapunov Analysis of Accelerated PDHG Algorithms}
\author[1]{Xueying Zeng}
\author[2,3]{Bin Shi\thanks{Corresponding author: \url{shibin@lsec.cc.ac.cn} } }
\affil[1]{School of Mathematical Sciences, Ocean University of China, Qingdao 266100, China}
\affil[2]{Academy of Mathematics and Systems Science, Chinese Academy of Sciences, Beijing 100190, China}
\affil[3]{School of Mathematical Sciences, University of Chinese Academy of Sciences, Beijing 100049, China}
\date\today

\maketitle

%\bin{Vanishing Viscosity should be included in the main part, the view is key important. At least, it is mentioned a little.}\wjs{OK.}
\begin{abstract}
The generalized \texttt{Lasso} is a remarkably versatile and extensively utilized model across a broad spectrum of domains, including statistics, machine learning, and image science. Among the optimization techniques employed to address the challenges posed by this model, saddle-point methods stand out for their effectiveness. In particular, the \textit{primal-dual hybrid gradient} (\texttt{PDHG}) algorithm has emerged as a highly popular choice, celebrated for its robustness and efficiency in finding optimal solutions. Recently,  the underlying mechanism of the~\texttt{PDHG} algorithm has been elucidated through the high-resolution ordinary differential equation (ODE) and the implicit-Euler scheme as detailed in~\citep{li2024understanding}. This insight has spurred the development of several accelerated variants of the ~\texttt{PDHG} algorithm, originally proposed by~\citep{chambolle2011first}. By employing discrete Lyapunov analysis, we establish that the~\texttt{PDHG} algorithm with iteration-varying step sizes, converges at a rate near $O(1/k^2)$. Furthermore, for the specific setting where $\tau_{k+1}\sigma_k = s^2$ and $\theta_k = \tau_{k+1}/\tau_k \in (0, 1)$ as proposed in~\citep{chambolle2011first}, an even faster convergence rate of $O(1/k^2)$ can be achieved. To substantiate these findings, we design a novel discrete Lyapunov function. This function is distinguished by its succinctness and straightforwardness, providing a clear and elegant proof of the enhanced convergence properties of the \texttt{PDHG} algorithm under the specified conditions. Finally, we utilize the discrete Lyapunov function to establish the optimal linear convergence rate when both the objective functions are strongly convex. 
\end{abstract}

%\thispagestyle{empty}
%\setcounter{page}{0}
%
%\newpage
%\tableofcontents
%\newpage

\section{Introduction}
\label{sec: intro}

In the realms of modern statistics and inverse problems, a widely embraced and highly effective approach is the incorporation of regularization penalties into model fitting~\citep{hoerl1970ridge}.  Among these methods, the $\ell_1$-penalty stands out as particularly significant, renowned for its prowess in capturing and promoting sparse structures.  This approach is epitomized by the~\textit{least absolute shrinkage and selection operator} (\texttt{Lasso})~\citep{tibshirani1996regression}, which also accounts for noise in the observed signal. Given $A \in \mathbb{R}^{m \times d_1}$, an $m \times d_1$ matrix, and $b \in \mathbb{R}^m$, a vector of  $m$ dimensions,  the~\texttt{Lasso} formulation is expressed as:
\[
\min_{x \in \mathbb{R}^{d_1}} \Phi(x) : = \frac{1}{2}\|Ax - b\|^2 + \lambda\|x\|_1,
\]
where the regularization parameter $\lambda > 0$ balances the tradeoff between fidelity to the measurements and sensitivity to noise.\footnote{Throughout this paper, the notation $\|\cdot\|$ specifically refers to the $\ell_2$-norm or Euclidean norm, $\|\cdot\|_2$, and the $\ell_1$-norm  is defined as
\[ \|x\|_1 = \sum_{i=1}^{d} |x_i|,\]
for any $x \in \mathbb{R}^d$. It is worth noting that the subscript $2$ is often omitted unless otherwise noted.} For a broader range of applications,  the
commonly used form is the generalized \texttt{Lasso}~\citep{tibshirani2011solution}, which is formulated as:
\[
\min_{x \in \mathbb{R}^{d_1}} \Phi(x) : = \frac{1}{2}\|Ax - b\|^2 + \lambda \|Fx\|_1,
\]
with $F \in \mathbb{R}^{d_2 \times d_1}$ characterizing the sparsity of $x$ in a transform domain.  This generalization significantly expands its applicability across various domains in machine learning and imaging science~\citep{tibshirani2011solution, chambolle2016introduction}, particularly supporting techniques such as~\textit{total-variation denoising}~\citep{rudin1992nonlinear}, \textit{fused lasso}~\citep{tibshirani2005sparsity}, $\ell_1$-\textit{trend filtering}~\citep{kim2009ell_1}, \textit{wavelet smoothing}~\citep{donoho1995adapting}, among others. The versatility and importance of the generalized \texttt{Lasso} are thus underscored in contemporary data analysis and modeling practices.

From an optimization standpoint, the generalized~\texttt{Lasso} can be elegantly formulated as a general optimization problem: 
\begin{equation}
\label{eqn: optimization}
\min_{x \in \mathbb{R}^{d_1}} \Phi(x) := f(x) + g(Fx).
\end{equation}
By employing the conjugate transformation, also known as the Legendre-Fenchel transformation, we can reformulate the optimization problem~\eqref{eqn: optimization} into the following minimax form:
\begin{equation}
\label{eqn: saddle-opt}
\min_{y \in \mathbb{R}^{d_2}} \max_{x \in \mathbb{R}^{d_1}} \Phi(x, y) = \max_{x \in \mathbb{R}^{d_1}} \min_{y \in \mathbb{R}^{d_2}} \Phi(x, y) := f(x) + \langle Fx, y \rangle - g^{\star}(y),
\end{equation}
a formulation pioneered by~\citet{arrow1958studies}. The minimax form~\eqref{eqn: saddle-opt} allows for the practical implementation of proximal operations. By incorporating a momentum step into the two proximal operations, we arrive at the enhanced \textit{primal-dual hybrid gradient} (\texttt{PDHG}) algorithm, described as follows: 
\begin{subequations}
\begin{empheq}[left=\empheqlbrace]{align}
  & x_{k+1} = \argmin_{x \in \mathbb{R}^{d_1}} \left\{ f(x) + \big\langle Fx, y_k \big\rangle + \frac{1}{2\tau_k} \|x - x_k \|^2 \right\},                                                  \label{eqn: pdhg-descent} \\
  & \overline{x}_{k+1} = x_{k+1} + \theta_k (x_{k+1} - x_k),                                                                                                                                                                        \label{eqn: pdhg-special} \\
  & y_{k+1} = \argmax_{y \in \mathbb{R}^{d_2}} \left\{ - g^{\star}(y) + \big \langle F\overline{x}_{k+1}, y \big\rangle - \frac{1}{2\sigma_k} \|y - y_k \|^2 \right\},      \label{eqn: pdhg-ascent}
\end{empheq}
\end{subequations}
as initially proposed by~\citet{chambolle2011first}. Here, $\tau_k$ and $\sigma_k$ are step sizes that are dependent on iterations, and $\theta_k \in (0, 1]$ is a newly added parameter also varying with iterations. This adaptability potentially leads to faster convergence rates and improved performance in various practical applications. The \texttt{PDHG} algorithm has been widely employed to solve optimization problem~\eqref{eqn: optimization} when the proximal operations \eqref{eqn: pdhg-descent} and \eqref{eqn: pdhg-ascent} can be computed efficiently.

 In the optimization problem~\eqref{eqn: optimization}, it is common to assume that both functions, $f$ and $g$, are convex.  Accordingly, the parameters in the~\texttt{PDHG} algorithm are set to fixed values: $\tau_k = \tau$, $\sigma_k = \sigma$, $\tau\sigma = s^2$, and $\theta_k=1$. The foundational convergence analysis, originally proposed by~\citet{chambolle2011first}, has since been refined and expanded in subsequent studies~\citep{he2012convergence, he2014convergence, he2022generalized}.  Recently, the high-resolution ordinary differential equation (ODE) framework has been introduced to provide a deeper insight into Nesterov's acceleration~\citep{shi2022understanding} and the \texttt{ADMM} algorithm~\citep{li2024understanding1}. This innovative approach, coupled with the implicit-Euler scheme, has been successfully elucidated the underlying mechanisms of~\texttt{PDHG} in the aforementioned scenario~\citep{li2024understanding}.  This advancement marks a significant leap in understanding the complexities of optimization algorithms and their convergence properties, offering profound insights into their dynamic behavior and performance characteristics. In this paper, we expand on the use of Lyapunov analysis to understand and discuss the convergence rates of accelerated~\texttt{PDHG} algorithms with iteration-varying parameters and various objective functions.  By exploring these variations, we aim to provide a comprehensive understanding of how these parameters influence the convergence rates and overall efficiency of the algorithms, thereby contributing to the development of more robust and versatile optimization techniques.

%from fixed $\tau_k = \tau$, $\sigma_k = \sigma$, $\tau\sigma = s^2$ and $\theta_k=1$ to iteration-varying parameters  $\tau_k \sigma_k = s^2$ and $\theta_k = 1$

% In this process, it is observed that the time $t_k = ks$ actually has only one scale $s$. 

%%%%%%%%%%%%%%%%%%%%%%%%%%%%%%%%%%%%%%%%%%%%%%%%%%%%%%%%%%%%%%%%%%%%%%%%%%%%%%%%%%%%%%%%%%%%%%%%%%%%%%%%%%%%%%%%%%%
\subsection{One advantage of discrete Lyapunov functions}
\label{subsec: merit-discrete-lyapunov}

To fully appreciate the advantages of discrete Lyapunov functions, we first assume that both objective functions, $f$ and $g$, are sufficiently smooth. Consider the scenario where $\tau_k = \tau$, $\sigma_k = \sigma$, $\tau\sigma = s^2$ and $\theta_k=1$. We elaborate on the~\texttt{PDHG} algorithm,~\eqref{eqn: pdhg-descent} -~\eqref{eqn: pdhg-ascent}, as follows:
\begin{subequations}
\begin{empheq}[left=\empheqlbrace]{align}
  & \frac{x_{k+1} - x_{k}}{\tau}      - F^{\top} (y_{k+1} - y_{k}) = - F^{\top} y_{k+1} - \nabla f(x_{k+1}),        \label{eqn: pdhg1-descent} \\ 
  & \frac{y_{k+1} - y_{k}}{\sigma} - F (x_{k+1} - x_{k})             =   Fx_{k+1}   - \nabla g^{\star}(y_{k+1}).      \label{eqn: pdhg1-ascent}
\end{empheq}
\end{subequations}
Reflecting on the process of deriving the high-resolution ODE, it becomes necessary to introduce a new parameter $s$ as the step size. By utilizing the $O(s)$-order approximation, we derive the system of high-resolution ODEs:
\begin{subequations}
\begin{empheq}[left=\empheqlbrace]{align}
  & \frac{s}{\tau}\dot{X} - s F^{\top} \dot{Y} = - F^{\top}Y - \nabla f(X),                          \label{eqn: high-ex-descent} \\
  & \frac{s}{\sigma}\dot{Y} - s F\dot{X} =  FX - \nabla g^{\star}(Y).                           \label{eqn: high-ex-ascent}
\end{empheq}
\end{subequations}
The continuous Lyapunov function, as formulated in~\citep[(6.3)]{li2024understanding}, is expressed as 
\begin{equation}
\label{eqn: lyapunov-ex-ode}
\mathcal{E}(t) = \frac{1}{2\tau} \|X - x^{\star}\|^2 + \frac{1}{2\sigma} \|Y - y^{\star}\|^2 - \big\langle F(X - x^{\star}), Y - y^{\star} \big\rangle. 
\end{equation}
By calculating the time derivative $\frac{d\mathcal{E}}{dt}$, we determine the convergence rate. Notably, the time $t_k = ks$ exhibits a single scale $s$. In contrast, for the discrete~\texttt{PDHG} algorithm,~\eqref{eqn: pdhg1-descent} and~\eqref{eqn: pdhg1-ascent},  we encounter dual scales, $\tau$ and $\sigma$. The discrete Lyapunov function, introduced in~\citep[(6.4)]{li2024understanding}, is formulated as:
\begin{equation}
\label{eqn: pdhg-ex-lyapunov}
\mathcal{E}(k) = \frac{1}{2\tau} \|x_k - x^{\star}\|^2 + \frac{1}{2\sigma} \|y_k - y^{\star}\|^2 - \big\langle F(x_k - x^{\star}), y_k - y^{\star} \big\rangle.
\end{equation}
This formulation facilitates the direct computation of the iterative difference $\mathcal{E}(k+1) - \mathcal{E}(k)$.  Thus, while the time derivative of the continuous Lyapunov function pertains a single scale, the iterative difference of the discrete Lyapunov function accommodate multiple scales. Consequently, in our analysis, we calculate the iterative difference of the discrete Lyapunov function, which can encompass multiple scales.  Building on this advantage, we proceed to employ discrete Lyapunov functions for analyzing iteration-varying parameters, thereby emphasizing their effectiveness in evaluating convergence rates.

%%%%%%%%%%%%%%%%%%%%%%%%%%%%%%%%%%%%%%%%%%%%%%%%%%%%%%%%%%%%%%%%%%%%%%%%%%%%%%%%%%%%%%%%%%%%%%%%%%%%%%%%%%%%%%%%%%%
\subsection{Acceleration: iteration-varying $\theta_k =\tau_{k+1}/\tau_{k}\in (0, 1)$ and $\tau_{k+1} \sigma_{k} = s^{2}$ }
\label{subsec: para_theta}

In practical applications, it is often more reasonable to assume the function $f$ in the composite optimization problem~\eqref{eqn: optimization} is strongly convex. This assumption is well-illustrated by models such as~\texttt{Lasso} and genralized~\texttt{Lasso}. By leveraging this strong convexity, the~\texttt{PDHG} algorithm can be significantly enhanced to achieve a faster convergence rate.  A promising approach involves considering the~\texttt{PDHG} algorithm with iteration-varying step sizes, $\tau_k$ and $\sigma_k$, while keeping the parameter $\theta_k = 1$ fixed:
\begin{subequations}
\begin{empheq}[left=\empheqlbrace]{align}
  & \frac{x_{k+1} - x_{k}}{\tau_k}      - F^{\top} (y_{k+1} - y_{k}) = - F^{\top} y_{k+1} - \nabla f(x_{k+1}),          \label{eqn: pdhg1-vary-descent} \\ 
  & \frac{y_{k+1} - y_{k}}{\sigma_k} - F (x_{k+1} - x_{k})             =   Fx_{k+1}   - \nabla g^{\star}(y_{k+1}).        \label{eqn: pdhg1-vary-ascent}
\end{empheq}
\end{subequations}
In alignment with this approach, the Lyapunov function is also enhanced with iteration-varying parameters as:
\begin{equation}
\label{eqn: pdhg-ex-lyapunov-very}
\mathcal{E}(k) = \frac{1}{2\tau_k} \|x_k - x^{\star}\|^2 + \frac{1}{2\sigma_k} \|y_k - y^{\star}\|^2 - \big\langle F(x_k - x^{\star}), y_k - y^{\star} \big\rangle.
\end{equation}
This methodology indeed facilitates a faster convergence rate, though it does not quite achieve the faster rate of $O(1/k^2)$. However, the enhanced convergence rate brought about by iteration-varying parameters represents a significant improvement over fixed parameter approaches. Further details and in-depth analysis of this approach can be found in~\Cref{sec: strong-con-1}.

\citet[Section 5.1]{chambolle2011first} introduced an innovative setting for iteration-varying parameters, notably incorporating $\theta_k \in (0, 1)$ to finely tune the convergence dynamics of the~\texttt{PDHG} algorithm.  Their seminal work also involved refining the relationship between $\tau_k$ and $\sigma_k$ from $\tau_{k} \sigma_{k} = s^2$ to $\tau_{k+1} \sigma_{k} = s^2$. This strategic adjustment has proven instrumental in elevating the efficiency of the~\texttt{PDHG} algorithm, enabling it to achieve the faster convergence rate of $O(1/k^2)$. In this paper, we propose a novel discrete Lyapunov function designed to provide a concise proof of convergence:
\begin{equation}
\label{eqn: lyapunov-acceleration}
\mathcal{E}(k) = \frac{\|x_k - x^{\star}\|^2}{2\tau_k^2} + \frac{\|y_{k-1} - y^{\star}\|^2}{2s^2} + \frac{\left\langle F(x_k - x_{k-1}), y_{k-1} - y^{\star} \right\rangle}{\tau_{k-1}} + \frac{\|x_k - x_{k-1}\|^2}{2\tau_{k-1}^2}.
\end{equation}
This discrete Lyapunov function is meticulously crafted to elucidate the iterative progression and convergence properties of the algorithm. For a comprehensive understanding, the detailed proof is elaborated in~\Cref{sec: strong-con-2}, showcasing the efficacy and applicability of our proposed approach in optimization theory.

%%%%%%%%%%%%%%%%%%%%%%%%%%%%%%%%%%%%%%%%%%%%%%%%%%%%%%%%%%%%%%%%%%%%%%%%%%%%%%%%%%%%%%%%%%%%%%%%%%%%%%%%%%%%%%%%%%%
\subsection{Overview of contributions}
\label{subsec: overview}

In this paper, we develop and employ discrete Lyapunov functions to rigorously establish the convergence rates for accelerated~\texttt{PDHG} algorithms with iteration-varying parameters in a diverse range of objective functions. Our key contributions are highlighted as follows. 

\begin{itemize}
\item We begin by analyzing the accelerated~\texttt{PDHG} algorithm with iteration-varying step sizes, specifically focusing on the case where the objective function $f$ is strongly convex. To this end, we refine the discrete Lyapunov function~\eqref{eqn: pdhg-ex-lyapunov}, originally proposed in~\citep{li2024understanding}, to account for step sizes that vary with each iteration, as outlined in~\eqref{eqn: pdhg-ex-lyapunov-very}.  A reasonable choice for these step sizes is $\tau_k \propto (k+1)^{-1}$ and $\sigma_k \propto (k+1)$, with the parameter $\theta_k = 1$ fixed. This leads to a convergence rate:
\[
\|x_k - x^{\star}\|^2 \leq O\left( \frac{1}{k^{2-\epsilon}}\right),
\]
where $\epsilon >0$ can be arbitarily small.  

\item In the seminal work by~\citet{chambolle2011first},  the step sizes are set such that $\tau_{k+1}\sigma_{k} = s^2$ and the parameter $\theta_k = \tau_{k+1}/\tau_k \in (0,1)$ , resulting in the faster convergence rate:
\[
\|x_k - x^{\star}\|^2 \leq O\left( \frac{1}{k^2}\right).
\]
In this paper, we propose a novel discrete Lyapunov function as given in~\eqref{eqn: lyapunov-acceleration}. By selecting the step sizes $\tau_k \propto k^{-1}$ and $\sigma_k \propto k+1$, along with the parameter $\theta_{k} =\frac{k}{k+1}$, we present a proof for the accelerated convergence rate.  Our proof, compared to the earliear work~\citep{chambolle2011first}, is distinguished by its succinctness and straightforwardness.

\item Finally, we delve into the analysis of the~\texttt{PDHG} algorithm, concentrating on the scenarios where both objective functions, $f$ and $g^{\star}$, exhibit strong convexity.  By utilizing the discrete Lyapunov function~\eqref{eqn: pdhg-ex-lyapunov}, we establish linear convergence. Furthermore, we propose the optimal configuration of step sizes $\tau$ and $\sigma$ to achieve the best possible convergence rate.

%We employ the the discrete Lyapunov function~\eqref{eqn: pdhg-ex-lyapunov} to establish the linear convergence. Finally, we propose the step sizes settings $\tau$ and $\sigma$ to achieve  the optimal rate. 
\end{itemize}

\section{Preliminaries}
\label{sec: prelim}

In this section, we embark on a concise journey through basic definitions and classical theorems in convex optimization. These definitions and theorems are primarily drawn from the classical literature, such as~\citep{rockafellar1970convex, rockafellar2009variational, nesterov1998introductory, boyd2004convex}, and serve as crucial references for proofs and discussions in the subsequent sections. 

An extended real-valued function is a function that maps to the extended real line $[-\infty,\infty]$. Unless otherwise specified, the functions in this work are extended valued. For an extended valued function $f$, its effective domain is the set
$$\mathrm{dom}(f)=\{x\in\mathbb{R}^d|f(x)<\infty\}.$$
A function is proper if its value is never $-\infty$ and is finite somewhere.

Let us begin by defining a convex function, its conjugate function, and the concepts of subgradient and subdifferential. 
\begin{defn}%[Definition 3.1.1 in \citep{nesterov1998introductory}]
\label{defn: convex} 
A function $f: \mathbb{R}^{d} \mapsto (-\infty,\infty]$ is said to be convex if, for  any $x, y \in \mathbb{R}^d$ and any $\lambda \in [0,1]$, the following inequlaity holds:
\[
f\left( \lambda x + (1 - \lambda)y \right) \leq \lambda f(x) + (1 - \lambda)f(y).
\] 
The function $f^{\star}: \mathbb{R}^{d} \mapsto  (-\infty,\infty]$ is said to be the conjugate of $f$ if, for any $y \in \mathbb{R}^d$,  it satisfies the following definition:
\[
f^{\star}(y):= \sup_{x \in \mathbb{R}^{d}} \left( \langle y, x \rangle - f(x)\right). 
\]
\end{defn}

\begin{defn}%[Definition 3.1.6 in \citep{nesterov1998introductory}]%[\citep[Definition 3.1.6]{nesterov1998introductory}]
\label{defn: subgradient}
 A vector $v$ is said to be a subgradient of a convex function $f$ at $x\in\mathrm{dom}(f)$ if, for any $y \in \mathbb{R}^d$, the following inequality holds:
\[
f(y) \geq f(x) + \langle v, y - x\rangle.
\] 
Additionally, the collection of all subgradients of $f$ at $x$ is the subdifferential of $f$ at $x$ and is denoted as $\partial f(x)$. 
\end{defn}

We can establish a sufficient and necessary condition for a point to be a minimizer of any convex function. 

\begin{theorem}
\label{thm: minimizer} Let $f$ be a proper convex function. Then $x^{\star}$ is a minimizer of $f$ if and only if $0 \in \partial f(x^{\star})$. 
\end{theorem}

%Next, we outline a basic property of subdifferentials, as stated in~\citep{rockafellar1970convex}. 

%\begin{theorem}[Theorem 23.8 in \citep{rockafellar1970convex}]
%\label{thm: sum-diff}
%For any $f_1, f_2 \in\mathcal{F}(\mathbb{R}^d)$, the subdifferentials satisfy
%\[
%\partial(f_1 + f_2)(x) = \partial f_1(x) + \partial f_2(x), 
%\]
%for any $x \in \mathbb{R}^d$. 
%\end{theorem}

%\begin{theorem}[Theorem 25.1 in \citep{rockafellar1970convex}]
%\label{thm: subgrad-unique}
%For any $f \in \mathcal{F}^1(\mathbb{R}^d)$, the gradient $\nabla f(x)$ is the unique subgradient of $f$ at $x$, satisfying 
%\[
%f(y) \geq f(x) + \langle \nabla f(x), y - x\rangle,
%\] 
%for any $y \in \mathbb{R}^d$.
%\end{theorem}

Let $\mathcal{F}(\mathbb{R}^d)$ be the class of proper, lower semi-continuous, and convex functions. For any function in $\mathcal{F}(\mathbb{R}^d)$, its conjugate belongs to $\mathcal{F}(\mathbb{R}^d)$. We denote $\mathcal{S}_{\mu}(\mathbb{R}^d) \subset \mathcal{F}(\mathbb{R}^d)$ as the class of $\mu$-strongly convex functions.

%For any $f \in \mathcal{F}(\mathbb{R}^d)$, we can provide a precise characterization as follows. 
% Next, we define the $\mu$-strongly convex function instead of the general convex function for the smooth part in the composite function to characterize the least-square model. 
\begin{defn}%[Definition 3.1.1 in \citep{nesterov1998introductory}]
\label{defn: strongly}
Let $f \in \mathcal{F}(\mathbb{R}^d)$ and $\mu>0$. For any $x, y \in \mathbb{R}^d$, if the function $f$ satisfies
\[
f(y) \geq f(x) + \langle v, y - x \rangle + \frac{\mu}{2} \|y - x\|^2, 
\]
where $v \in \partial f(x)$ is a subgradient, we say that $f$ is $\mu$-strongly convex. 
\end{defn}
%\begin{itemize}
%\item[(1)] if the gradient of the function $f$ satisfies
%\[
%\|\nabla f(x) - \nabla f(y)\| \leq L \|x - y\|,
%\]
%we say that $f$ is $L$-smooth;
%\item[(2)] if the function $f$ satisfies 
%\[
%f(y) \geq f(x) + \langle \nabla f(y), y - x \rangle + \frac{\mu}{2} \|y - x\|^2, 
%\]
%we say that $f$ is $\mu$-strongly convex. 
%\end{itemize}

%Given that the implicit solution in the iteration~\eqref{eqn: pdhg-descent} may not always be directly obtainable,  especially in scenarios where $f$ is non-quadratic, it is necessary to specify the class of the objective function $f$ that allows for the direct obtainment of the implicit solution in these iterations. Let us denote this class as $\mathcal{R}(\mathbb{R}^d)$.  
Finally, we introduce the definition of saddle points as denoted in~\citep[Definition 11.49]{rockafellar2009variational}. 
\begin{defn}%[Definition 3.1.6 in \citep{nesterov1998introductory}]%[\citep[Definition 3.1.6]{nesterov1998introductory}]
\label{defn: saddle} Let $\Phi:\mathbb{R}^{d_1}\times\mathbb{R}^{d_2}\rightarrow[-\infty,+\infty]$.
 A vector pair $(x^{\star},y^{\star}) \in \mathbb{R}^{d_1} \times \mathbb{R}^{d_2}$ is said to be a saddle point of the convex-concave function $\Phi$ (convex with respect to the variable $x$ and concave with respect to the variable $y$) if, for any $x \in \mathbb{R}^{d_1}$ and $y \in \mathbb{R}^{d_2}$, the following inequality holds: 
\[
\Phi(x^{\star},y) \leq \Phi(x^{\star},y^{\star}) \leq  \Phi(x,y^{\star}).
\]
\end{defn}
With~\Cref{defn: saddle},  it becomes straightforward for us to derive the sufficient and necessary condition for the saddle point of the objective function $\Phi(x,y)$ taking the convex-concave form~\eqref{eqn: saddle-opt}. We conclude this section with the following theorem.
\begin{theorem}
\label{thm: saddle-sn}
Let $f \in \mathcal{F}(\mathbb{R}^{d_1})$, $g \in \mathcal{F}(\mathbb{R}^{d_2})$ and $F \in \mathbb{R}^{d_2 \times d_1}$ be a $d_2 \times d_1$ matrix. Given that the objective function $\Phi(x,y)$ takes the convex-concave form~\eqref{eqn: saddle-opt},  the point $(x^{\star},y^{\star})$ is a saddle point if and only if, for any $x \in \mathbb{R}^{d_1}$ and $y \in \mathbb{R}^{d_2}$, the following inequalities hold:
%\begin{subequations}
%\begin{empheq}[left=\empheqlbrace]{align}
\[
\left\{\begin{aligned}
          & f(x) - f(x^{\star}) + \big\langle F^\top y^\star,x - x^{\star} \big\rangle \geq 0,                            \\ %  \label{eqn: saddle-x} \\
          & g^{\star}(y) - g^{\star}(y^{\star}) -  \big\langle  Fx^{\star}, y - y^{\star} \big\rangle \geq 0.         % \label{eqn: saddle-y}   
         \end{aligned} \right.
\]        
Furthermore, the functions, $f$ and $g^{\star}$, at the saddle point $(x^{\star}, y^{\star})$, satisfy the following conditions:
\[
\left\{\begin{aligned}
          &   \partial f(x^{\star})              +  F^{\top}y^{\star} \ni 0,                            \\ %  \label{eqn: saddle-x} \\
          &   \partial g^{\star}(y^{\star}) -   Fx^{\star}            \ni 0.         % \label{eqn: saddle-y}   
         \end{aligned} \right.
\]          
%\end{empheq}
%\end{subequations}
\end{theorem}

\section{The case: $f \in \mathcal{S}_{\mu}(\mathbb{R}^{d_1})$}
\label{sec: strong-con-1}

%\begin{equation}
%\label{eqn: lyapunov}
%\mathcal{E}(k) = \frac{ \| x_k - x^{\star} \|^2 }{2\tau_k}+ \frac{\|y_k - y ^{\star}\|^2}{2\sigma_k}  - \left\langle F(x_k - x^{\star}), y_k - y^{\star} \right\rangle
%\end{equation}

In this section, we delve into the convergence rate for the~\texttt{PDHG} algorithm, given by~\eqref{eqn: pdhg-descent} -~\eqref{eqn: pdhg-ascent}, utilizing the iteration-varying discrete Lyapunov function~\eqref{eqn: pdhg-ex-lyapunov-very}. This analysis considers iteration-varying step sizes, $\tau_k$ and $\sigma_k$, with a fixed parameter $\theta_k = 1$. By computing the iterative difference, we obtain the following: 
\begin{align}
& \mathcal{E}(k+1)  - \mathcal{E}(k)  \nonumber  \\
& = \left\langle \frac{x_{k+1} - x_{k}}{\tau_k} - F^{\top}(y_{k+1} - y_{k}), x_{k+1} - x^{\star} \right\rangle  + \left\langle \frac{y_{k+1} - y_{k}}{\sigma_k} - F(x_{k+1} - x_{k}), y_{k+1} - y^{\star} \right\rangle \nonumber \\
& \mathrel{\phantom{=}} + \left( \frac{1}{2\tau_{k+1}} - \frac{1}{2\tau_k} \right) \|x_{k+1} - x^{\star}\|^2 + \left( \frac{1}{2\sigma_{k+1}} - \frac{1}{2\sigma_k} \right) \|y_{k+1} - y^{\star}\|^2                     \nonumber \\
& \mathrel{\phantom{=}} - \underbrace{\left( \frac{1}{2\tau_k} \|x_{k+1} - x_{k}\|^2 + \frac{1}{2\sigma_k} \|y_{k+1} - y_{k}\|^2 - \left\langle F(x_{k+1} - x_{k}), y_{k+1} - y_{k} \right\rangle \right)}_{\mathbf{NE}}, \label{eqn: iter-diff-1}
\end{align}
where $\mathbf{NE}$ represents the numerical error resulting from the implicit discretization.  With the properties of objective functions and step sizes, we can further estimate the iterative difference~\eqref{eqn: iter-diff-1} in the following lemma.

\begin{lemma}
\label{lem: strong-con-1}
Let $f \in \mathcal{S}_{\mu}(\mathbb{R}^{d_1})$ and $g \in \mathcal{F}(\mathbb{R}^{d_2})$.  For any step size $\tau_k \sigma_k = s^2 \in (0, \|F\|^{-2})$ ($k=0,1,2,\ldots,n,\ldots$), the discrete Lyapunov function~\eqref{eqn: pdhg-ex-lyapunov-very} satisfies the following inequality:
\begin{equation}
\label{eqn: strong-con-11}
\mathcal{E}(k+1) - \mathcal{E}(k) \leq - \left( \mu + \frac{1}{2\tau_k} - \frac{1}{2\tau_{k+1}} \right) \|x_{k+1} - x^{\star}\|^2 -  \left( \frac{1}{2\sigma_k} - \frac{1}{2\sigma_{k+1}} \right) \|y_{k+1} - y^{\star}\|^2.
\end{equation}
\end{lemma}

\begin{proof}[Proof of~\Cref{lem: strong-con-1}]
For each $k = 0, 1, 2, \ldots, n, \ldots$,  the relation $\tau_k \sigma_k = s^2 \in (0, \|F\|^{-2})$ holds. By applying the Cauchy-Schwarz inequality, we establish the numerical error $\mathbf{NE} \geq 0$. Hence, we can estimate the iterative difference~\eqref{eqn: iter-diff-1} as:
\begin{align}
& \mathcal{E}(k+1)  - \mathcal{E}(k)  \nonumber  \\
& \leq \left\langle \frac{x_{k+1} - x_{k}}{\tau_k} - F^{\top}(y_{k+1} - y_{k}), x_{k+1} - x^{\star} \right\rangle  + \left\langle \frac{y_{k+1} - y_{k}}{\sigma_k} - F(x_{k+1} - x_{k}), y_{k+1} - y^{\star} \right\rangle \nonumber \\
& \mathrel{\phantom{=}} + \left( \frac{1}{2\tau_{k+1}} - \frac{1}{2\tau_k} \right) \|x_{k+1} - x^{\star}\|^2 + \left( \frac{1}{2\sigma_{k+1}} - \frac{1}{2\sigma_k} \right) \|y_{k+1} - y^{\star}\|^2.        \label{eqn: iter-diff-2}
\end{align}
Given any $f \in \mathcal{S}_{\mu}(\mathbb{R}^{d_1})$ and any $g \in \mathcal{F}(\mathbb{R}^{d_1})$,~\Cref{thm: minimizer} allows us to bypass the argmin and argmax operations, enabling further exploration of the iteration-varying~\texttt{PDHG} algorithm,~\eqref{eqn: pdhg-descent} -~\eqref{eqn: pdhg-ascent}. Specifically, for the first descent iteration~\eqref{eqn: pdhg-descent} and the third ascent iteration~\eqref{eqn: pdhg-ascent}, we express: 
\begin{subequations}
\begin{empheq}[left=\empheqlbrace]{align}
  & \frac{x_{k+1} - x_{k}}{\tau_k}      - F^{\top} (y_{k+1} - y_{k}) + F^{\top} y_{k+1} + \partial f(x_{k+1}) \ni 0,        \label{eqn: pdhga1-descent} \\ 
  & \frac{y_{k+1} - y_{k}}{\sigma_k} - F (x_{k+1} - x_{k})   -  Fx_{k+1}   + \partial g^{\star}(y_{k+1}) \ni 0.                \label{eqn: pdhga1-ascent}
\end{empheq}
\end{subequations}
According to~\Cref{thm: saddle-sn},  the saddle point $(x^{\star}, y^{\star})$ of the objective function $\Phi(x, y)$ in the form~\eqref{eqn: saddle-opt} satisfies:
\begin{subequations}
\begin{empheq}[left=\empheqlbrace]{align}
          &   \partial f(x^{\star})              +  F^{\top}y^{\star} \ni 0,            \label{eqn: saddle-x-defn} \\
          &   \partial g^{\star}(y^{\star}) -   Fx^{\star}            \ni 0.            \label{eqn: saddle-y-defn}   
\end{empheq}
\end{subequations}
By \Cref{defn: strongly} of strongly convex functions, we can drive 
\begin{equation*}
	f(x^\star)-f(x^{k+1})\geq \left<\frac{x_{k+1} - x_{k}}{\tau_k}      - F^{\top} (y_{k+1} - y_{k}) + F^{\top} y_{k+1},x_{k+1}-x^\star\right>+\frac{\mu}{2}\|x_{k+1}-x^\star\|^2,
\end{equation*}
from \eqref{eqn: pdhga1-descent} and
\begin{equation*}
	f(x^\star)-f(x^{k+1})\leq \left<F^\top y^\star,x_{k+1}-x^\star\right>-\frac{\mu}{2}\|x_{k+1}-x^\star\|^2,
\end{equation*}
from \eqref{eqn: saddle-x-defn}. Therefore, we deduce that
\begin{equation}\label{eqn:strongconcexeq}
 \left<\frac{x_{k+1} - x_{k}}{\tau_k}      - F^{\top} (y_{k+1} - y_{k}) + F^{\top} (y_{k+1}-y^\star),x_{k+1}-x^\star\right>\leq-\mu\|x_{k+1}-x^\star\|^2.
\end{equation}
Similarly, by \Cref{defn: subgradient} of subgradient, we can derive
$$g^\star(y^\star)-g^\star(y_{k+1})\geq\left\langle \frac{y_{k+1} - y_{k}}{\sigma_k} - F (x_{k+1} - x_{k})   -  Fx_{k+1} ,y_{k+1}-y^\star\right\rangle,$$
and
$$g^\star(y^\star)-g^\star(y_{k+1})\leq-\left\langle Fx^\star ,y_{k+1}-y^\star\right\rangle,$$
from
\eqref{eqn: pdhga1-ascent} and \eqref{eqn: saddle-y-defn}, respectively. By the last two inequalities, we deduce that
\begin{equation}\label{eqn:concexeq}
	 \left\langle \frac{y_{k+1} - y_{k}}{\sigma_k} - F(x_{k+1} - x_{k}) - F(x_{k+1} - x^{\star}), y_{k+1} - y^{\star} \right\rangle\leq 0.
\end{equation}
By summing \eqref{eqn:strongconcexeq} and \eqref{eqn:concexeq} and substituting the result into \eqref{eqn: iter-diff-2}, we complete the proof.
\end{proof}

Let us introduce a positive constant $c > 0$. Using the Cauchy-Schwarz inequality, we can estimate the discrete Lyapunov function~\eqref{eqn: pdhg-ex-lyapunov-very} as follows:
\begin{equation}
\label{eqn: estimate-lyapunov}
\mathcal{E}(k+1) \leq \frac{1}{2\tau_{k+1}}\left( 1 + \frac{s}{c} \right) \|x_{k+1} - x^{\star}\|^2 + \frac{1}{2\sigma_{k+1}} \left( 1 + cs\|F\|^2 \right) \|y_{k+1} - y^{\star}\|^2. 
\end{equation}
This inequality serves as an upper bound for $\mathcal{E}(k+1)$.  Consequently, by comparing it with the right side of  the inequality~\eqref{eqn: strong-con-11} and  scaling the corresponding terms, we can estimate the following proportion:
\[
\frac{\mathcal{E}(k+1) - \mathcal{E}(k)}{\mathcal{E}(k+1)}. 
\] 
Additionally, an intuitive observation suggests that the ratio of the coefficients of the terms $\|x_{k+1} - x^{\star}\|^2$ and $\|y_{k+1} - y^{\star}\|^2$ in~\eqref{eqn: strong-con-11} and~\eqref{eqn: estimate-lyapunov} likely scales $O(1/k)$. Hence, a reasonable choice for the step sizes is:
\begin{equation}
\label{eqn: tau-sigma}
\tau_{k} = \frac{1}{c(k+1)}, \qquad \sigma_{k} = cs^2 (k + 1).
\end{equation}
With the step sizes given in~\eqref{eqn: tau-sigma}, we achieve the following conclusion.
\begin{theorem}
\label{thm: strong-con-1}
Let $c \in (0, 2\mu)$ be a positive constant. By taking the step sizes as given in~\eqref{eqn: tau-sigma} and under the same assumptions as in~\Cref{lem: strong-con-1}, the discrete Lyapunov function~\eqref{eqn: pdhg-ex-lyapunov-very} holds for any $k \geq 0$ such that: 
\begin{equation}
\label{eqn: lyapunov-estimate-1}
\mathcal{E}(k) \leq \frac{(1+\alpha)(k+1)!}{(k + 1 + \alpha)!} \mathcal{E}(0),
\end{equation}
where $\alpha = \min \left( \frac{2\mu - c}{s + c}, \frac{1}{1+cs\|F\|^2}\right)$. Furthermore, the iterative sequence $\{x_{k}\}_{k=0}^{\infty}$ generated by the~\texttt{PDHG} algorithm,~\eqref{eqn: pdhg-descent} -~\eqref{eqn: pdhg-ascent}, converges to $x^{\star}$ at the rate characterized by
\begin{equation}
\label{eqn: final-estimate-1}
\|x_k - x^{\star}\|^2 \leq \frac{1+ s\|F\|}{1 - s\|F\|}  \cdot \frac{(1+\alpha) \cdot k!}{(k+1+\alpha)!}\left( \|x_0 - x^{\star}\|^2 + \frac{1}{c^2s^2}\|y_0 - y^{\star}\|^2 \right). 
\end{equation}
\end{theorem}

\begin{remark}
In~\Cref{thm: strong-con-1}, it is observed that if we set the step size such that $s < \min(2\mu, \|F\|^{-1})$, then $\alpha \rightarrow 1^{-}$ as $c \rightarrow 0^{+}$. Consequently, the convergence rate described in~\eqref{eqn: final-estimate-1} becomes very close to
\[
\|x_k - x^{\star}\|^2 \leq O\left( \frac{1}{k^2} \right). 
\]
\end{remark}

\begin{proof}[Proof of~\Cref{thm: strong-con-1}]
By substituting the step sizes given in~\eqref{eqn: tau-sigma} into the iterative difference~\eqref{eqn: strong-con-11} and the estimate inequality~\eqref{eqn: estimate-lyapunov}, we derive the following relations:
\begin{subequations}
\begin{empheq}[left=\empheqlbrace]{align}
\mathcal{E}(k+1) - \mathcal{E}(k)&  \leq - \left( \mu - \frac{c}{2} \right) \|x_{k+1} - x^{\star}\|^2 -  \frac{1}{2cs^2(k+1)(k+2)}\cdot \|y_{k+1} - y^{\star}\|^2,   \label{eqn: strong-con-12} \\
\mathcal{E}(k+1) & \leq \frac{(s+c)(k+2)}{2} \|x_{k+1} - x^{\star}\|^2 + \frac{(1 + cs\|F\|^2)}{2cs^2(k+2)}  \|y_{k+1} - y^{\star}\|^2.                                         \label{eqn: strong-con-13}
\end{empheq}
\end{subequations}
Taking the ratio of~\eqref{eqn: strong-con-12} and~\eqref{eqn: strong-con-13} and simplifying further, we get:

\begin{align}
\frac{\mathcal{E}(k+1) - \mathcal{E}(k)}{\mathcal{E}(k+1)} & \leq - \frac{1}{k+2} \cdot \frac{(2\mu - c) \|x_{k+1} - x^{\star}\|^2 + \frac{1}{cs^2(k+1)(k+2)} \|y_{k+1} - y^{\star}\|^2 }{(s + c) \|x_{k+1} - x^{\star}\|^2 + \frac{1+cs\|F\|^2}{cs^2(k+2)^2} \|y_{k+1} - y^{\star}\|^2}  \nonumber\\
                                                                                              & \leq - \frac{1}{k+2} \cdot \frac{(2\mu - c) \|x_{k+1} - x^{\star}\|^2 + \frac{1}{cs^2(k+2)^2} \|y_{k+1} - y^{\star}\|^2 }{(s + c) \|x_{k+1} - x^{\star}\|^2 + \frac{1+cs\|F\|^2}{cs^2(k+2)^2} \|y_{k+1} - y^{\star}\|^2}  \nonumber \\
                                                                                              & \leq - \frac{1}{k+2} \cdot  \min \left( \frac{2\mu - c}{s + c}, \frac{1}{1+cs\|F\|^2}\right). \label{eqn: strong-con-14}  
\end{align}
Taking $\alpha =  \min \left( \frac{2\mu - c}{s + c}, \frac{1}{1+cs\|F\|^2}\right)$, it is derived for the iterative inequality of the discrete Lyapunov function~\eqref{eqn: pdhg-ex-lyapunov-very} as
\[
\mathcal{E}(k+1) \leq \frac{k + 2}{k + 2 + \alpha} \mathcal{E}(k).
\]
Hence, we derive the estimate inequality for $\mathcal{E}(k)$ given in~\eqref{eqn: lyapunov-estimate-1}.  Furthermore, applying the Cauchy-Schwarz inequality, we can estimate the discrete Lyapunov function $\mathcal{E}(0)$ and $\mathcal{E}(k)$ given in~\eqref{eqn: pdhg-ex-lyapunov-very} as
\begin{equation}
\label{eqn: strong-con-15}
\mathcal{E}(0) \leq \frac{1+s\|F\|}{2}  \left(c \|x_{0} - x^{\star}\|^2 + \frac{1}{cs^2}\|y_{0} - y^{\star}\|^2 \right),
\end{equation}
and
\begin{align}
\label{eqn: strong-con-16}
\mathcal{E}(k) & \geq \frac{1-s\|F\|}{2} \left(c(k+1) \|x_{k} - x^{\star}\|^2 + \frac{1}{cs^2(k+1)}\|y_{k} - y^{\star}\|^2 \right) \nonumber \\
                        & \geq \frac{c(1-s\|F\| )(k+1) \|x_{k} - x^{\star}\|^2}{2},
\end{align}
respectively. By substituting~\eqref{eqn: strong-con-15} and~\eqref{eqn: strong-con-16} into the inequality~\eqref{eqn: lyapunov-estimate-1}, we complete the proof. 

\end{proof}

\section{Acceleration settings: iteration-varying $\theta_k = \tau_{k+1}/\tau_{k} \in (0, 1)$ and $\tau_{k+1} \sigma_{k} = s^{2}$}
\label{sec: strong-con-2}

In this section, we utilize the novel discrete Lyapunov function given in~\eqref{eqn: lyapunov-acceleration} to establish the accelerated convergence rate for the~\texttt{PDHG} algorithm,~\eqref{eqn: pdhg-descent} -~\eqref{eqn: pdhg-ascent}, with iteration-varying step sizes $\tau_{k+1}\sigma_k = s^2$ and an iteration-varying parameter $\theta_k \in (0, 1)$. By computing the iterative difference, we obtain:
\begin{align}
& \mathcal{E}(k+1) - \mathcal{E}(k) \nonumber \\
& = \frac{1}{\tau_k} \left\langle \frac{x_{k+1} - x_{k}}{\tau_k}, x_{k+1} - x^{\star} \right\rangle + \left( \frac{1}{2\tau_{k+1}^2} - \frac{1}{2\tau_k^2}  \right) \|x_{k+1} - x^{\star}\|^2  \nonumber \\
& \mathrel{\phantom{=}} + \frac{1}{\tau_k} \left\langle \frac{y_k - y_{k-1}}{\sigma_{k-1}}, y_{k} - y^{\star}  \right\rangle  - \frac{\langle F(x_{k} - x_{k-1}), y_{k} - y^{\star} \rangle}{\tau_{k-1}} + \frac{\langle F(x_{k+1} - x_{k}), y_{k} - y^{\star} \rangle}{\tau_k}  \nonumber \\ 
& \mathrel{\phantom{=}} - \underbrace{ \left( \frac{\|x_k - x_{k-1}\|^2}{2\tau_{k-1}^2} -  \frac{\langle F(x_{k} - x_{k-1}), y_{k} - y_{k-1} \rangle}{\tau_{k-1}}  + \frac{\|y_k - y_{k-1}\|^2}{2s^2} \right) }_{\mathbf{NE}},  \label{eqn: iter-diff-strong-con-2}
\end{align}
where $\mathbf{NE}$ represents the numerical error resulting from the implicit discretization. Utilizing the properties of objective function and the iteration-varying step size and parameters, we can further estimate the iterative difference~\eqref{eqn: iter-diff-strong-con-2} in the following lemma. 

\begin{lemma}
\label{lem: strong-con-2}
Let $f \in \mathcal{S}_{\mu}(\mathbb{R}^{d_1})$ and $g \in \mathcal{F}(\mathbb{R}^{d_2})$.  For any step sizes that satisfy $\tau_{k+1} \sigma_k = s^2 \in (0, \|F\|^{-2}
)$ and $\tau_{k+1} = \theta_{k} \tau_k$ ($k=0,1,2,\ldots,n,\ldots$), the discrete Lyapunov function~\eqref{eqn: lyapunov-acceleration} satisfies the following inequality:
\begin{equation}
\label{eqn: strong-con-21}
\mathcal{E}(k+1) - \mathcal{E}(k) \leq - \left( \frac{\mu}{\tau_k} + \frac{1}{2\tau_k^2} - \frac{1}{2\tau_{k+1}^2} \right) \|x_{k+1} - x^{\star}\|^2.
\end{equation}
\end{lemma}
\begin{proof}[Proof of~\Cref{lem: strong-con-2}]
For each $k = 0, 1, 2, \ldots, n, \ldots$,  the relation $\tau_{k+1} \sigma_k = s^2 \in (0, \|F\|^{-2})$ holds. Applying the Cauchy-Schwarz inequality, we establish the numerical error $\mathbf{NE} \geq 0$.  Taking the relationship $\tau_{k} = \theta_{k-1}\tau_{k-1}$, we can estimate the iterative difference~\eqref{eqn: iter-diff-strong-con-2} as follows:
\begin{align}
& \mathcal{E}(k+1)  - \mathcal{E}(k) \nonumber \\
& \leq \left( \frac{1}{2\tau_{k+1}^2} - \frac{1}{2\tau_k^2}  \right) \|x_{k+1} - x^{\star}\|^2 + \frac{1}{\tau_k} \left\langle \frac{x_{k+1} - x_{k}}{\tau_k}, x_{k+1} - x^{\star} \right\rangle \nonumber \\
& \mathrel{\phantom{=}} + \frac{1}{\tau_k} \left\langle \frac{y_k - y_{k-1}}{\sigma_{k-1}} - \theta_{k-1} F (x_{k} - x_{k-1}) , y_{k} - y^{\star}  \right\rangle  + \frac{\langle F(x_{k+1} - x_{k}), y_{k} - y^{\star} \rangle}{\tau_k}. \label{eqn: iter-diff-strong-22} 
\end{align}
By the optimal conditions of the first descent iteration \eqref{eqn: pdhg-descent} and the third ascent iteration \eqref{eqn: pdhg-ascent}, we have
%
%Given any $f \in \mathcal{S}_{\mu}(\mathbb{R}^{d_1}) \cap \mathcal{R}(\mathbb{R}^{d_1})$ and any $g \in \mathcal{F}(\mathbb{R}^{d_1})$,~\Cref{thm: minimizer} and~\Cref{thm: sum-diff} allow us to bypass the argmin and argmax operations, enabling further exploration of the iteration-varying~\texttt{PDHG} algorithm,~\eqref{eqn: pdhg-descent} ---~\eqref{eqn: pdhg-ascent}. Specifically, in terms of the first descent iteration~\eqref{eqn: pdhg-descent} and the third ascent iteration~\eqref{eqn: pdhg-ascent}, expressed as: 
\begin{subequations}
\begin{empheq}[left=\empheqlbrace]{align}
  & \frac{x_{k+1} - x_{k}}{\tau_k}      + F^{\top}  y_{k}  + \partial f(x_{k+1}) \ni 0,        \label{eqn: pdhga2-descent} \\ 
  & \frac{y_{k+1} - y_{k}}{\sigma_k} - \theta_k F (x_{k+1} - x_{k})   -  Fx_{k+1}   + \partial g^{\star}(y_{k+1}) \ni 0.                \label{eqn: pdhga2-ascent}
\end{empheq}
\end{subequations}
Since $f$ is $\mu$-strongly convex and $g^\star$ is convex, following the similar procedure in the proof of Lemma \ref{lem: strong-con-1}, we can derive the following inequalities:
%By~\Cref{defn: subgradient} and~\Cref{defn: strongly} of convex and strongly convex functions, the expansion of the~\texttt{PDHG} algorithm,~\eqref{eqn: pdhga2-descent} and \eqref{eqn: pdhga2-ascent}, and the relationship of saddle points,~\eqref{eqn: saddle-x-defn} and \eqref{eqn: saddle-y-defn}, lead to: 
\begin{subequations}
\begin{empheq}[left=\empheqlbrace]{align}
  & \left\langle \frac{x_{k+1} - x_{k}}{\tau_k} , x_{k+1} - x^{\star} \right\rangle \leq -\mu \|x_{k+1} - x^{\star}\|^2-\left\langle F^{\top}  (y_{k} -y^{\star}) ,x_{k+1} - x^{\star}\right\rangle,        \label{eqn: pdhga2-inq1} \\ 
  & \left\langle \frac{y_k - y_{k-1}}{\sigma_{k-1}} - \theta_{k-1} F (x_{k} - x_{k-1}) , y_{k} - y^{\star}  \right\rangle \leq \left\langle F(x_{k}-x^{\star}),y_k-y^\star\right\rangle.                          \label{eqn: pdhga2-inq2}
\end{empheq}
\end{subequations}
By substituting~\eqref{eqn: pdhga2-inq1} and~\eqref{eqn: pdhga2-inq2} into the iterative difference~\eqref{eqn: iter-diff-strong-22}, we complete the proof. 
\end{proof}

As we know, the step sizes satisfy $\tau_{k+1}\sigma_k = s^2 \in (0, \|F\|^{-2})$.  By utilizing the Cauchy-Schwarz inequality,  we can deduce that the discrete Lyapunov function~\eqref{eqn: lyapunov-acceleration} satisfies:
\begin{equation}
\label{eqn: strong-con-2-lyapunov-inq}
\mathcal{E}(k) \geq \frac{\|x_{k} - x^{\star}\|^2}{2\tau_{k}^2}. 
\end{equation}
If the discrete Lyapunov function $\mathcal{E}(k)$ is bounded above, the inequality~\eqref{eqn: strong-con-2-lyapunov-inq} leads to the following convergence rate:
\begin{equation}
\label{eqn: strong-con-2-lyapunov-inq1}
\|x_k - x^{\star}\|^2 \leq O(\tau_k^2). 
\end{equation}
From~\Cref{lem: strong-con-2}, we know that for the discrete Lyapunov function $\mathcal{E}(k)$ to remain bounded,  the following inequality must hold:
\begin{equation}
\label{eqn: strong-con-2-inq-coefficient}
\frac{\mu}{\tau_k} + \frac{1}{2\tau_k^2} - \frac{1}{2\tau_{k+1}^2} \geq 0. 
\end{equation}
To achieve the faster convergence rate $O(1/k^2)$, as indicated by~\eqref{eqn: strong-con-2-lyapunov-inq1},  it is reasonable to choose the step sizes as
\begin{equation}
\label{eqn: strong-con-2-step-size}
\tau_{k} = \frac{1}{ck}, \qquad \sigma_{k} = cs^{2}(k+1). 
\end{equation}
Substituting the setting of step sizes~\eqref{eqn: strong-con-2-step-size} into the inequality~\eqref{eqn: strong-con-2-inq-coefficient}, we can derive the condition:
\[
2\mu k- c \left( 2k + 1\right) \geq 0, 
\]
which simplifies to 
\[
k \geq \frac{c}{2\mu - 2c}. 
\]
Building on these deductions, we conclude this section with the following theorem. 
\begin{theorem}
\label{thm: strong-con-2}
Let $c \in (0, \mu)$ be a positive constant. By taking the step sizes as given in~\eqref{eqn: strong-con-2-step-size} and under the same assumption as in~\Cref{lem: strong-con-2}, there exists a positive constant $K_0 = \left \lceil \frac{c}{2\mu - 2c}  \right \rceil$ such that for any $k \geq K_0$, the iterative sequence $\{x_{k}\}_{k=K_0}^{\infty}$ generated by the~\texttt{PDHG} algorithm,~\eqref{eqn: pdhg-descent} -~\eqref{eqn: pdhg-ascent} converges with the following rate:
\begin{equation}
\label{sen: strong-con-thm-1}
\|x_{k} - x^{\star}\|^2 \leq \frac{2\mathcal{E}(K_0)}{c^2k^2}. 
\end{equation}
Specifically, if $c \leq \frac{2\mu}{3}$, then we have
\begin{equation}
\label{sen: strong-con-thm-2}
\|x_{k} - x^{\star}\|^2 \leq \frac{2}{c^2k^2} \left(  \frac{\|x_1 - x^{\star}\|^2}{2\tau_1^2} + \frac{\|y_{0} - y^{\star}\|^2}{2s^2} + \frac{\left\langle F(x_1 - x_{0}), y_{0} - y^{\star} \right\rangle}{\tau_{0}} + \frac{\|x_1 - x_{0}\|^2}{2\tau_{0}^2} \right). 
\end{equation}
\end{theorem}

%\[
% - \frac{\|x_k - x_{k-1}\|^2}{2\tau_{k-1}^2} +  \frac{\langle F(x_{k} - x_{k-1}), y_{k} - y_{k-1} \rangle}{\tau_{k-1}}  - \frac{\|y_k - y_{k-1}\|^2}{2s^2}  \leq 0
%\]
%
%$\tau_{k} = \tau_{k-1} \theta_{k-1}$

%\begin{subequations}
%\begin{empheq}[left=\empheqlbrace]{align}
%  & \frac{x_{k+1} - x_{k}}{\tau_k}     + F^{\top} y_{k} + \partial f(x_{k+1}) = 0,                                                                \label{eqn: pdhgb1-descent} \\ 
%  & \frac{y_{k+1} - y_{k}}{\sigma_k} - \theta_k F (x_{k+1} - x_{k})   -  Fx_{k+1}   + \partial g^{\star}(y_{k+1}) = 0.         \label{eqn: pdhgb1-ascent}
%\end{empheq}
%\end{subequations}

\section{The case: $f \in \mathcal{S}_{\mu}(\mathbb{R}^{d_1})$ and  $g^{\star} \in \mathcal{S}_{\gamma}(\mathbb{R}^{d_2})$}
\label{sec: strong-strong}

In this section, we leverage the discrete Lyapunov function provided in~\eqref{eqn: pdhg-ex-lyapunov} to establish the linear convergence of the~\texttt{PDHG} algorithm as described by~\eqref{eqn: pdhg-descent} -~\eqref{eqn: pdhg-ascent}. We focus on the scenario where the step sizes $\tau_k = \tau$ and $\sigma_k = \sigma$ are fixed, and the parameter $\theta_k = 1$ remains constant. Additionally, we consider the setting where $f \in \mathcal{S}_{\mu}(\mathbb{R}^{d_1})$ and  $g^{\star} \in \mathcal{S}_{\gamma}(\mathbb{R}^{d_2})$. By computing the iterative difference, we obtain the following expression:
\begin{align}
& \mathcal{E}(k+1)  - \mathcal{E}(k)  \nonumber  \\
& = \left\langle \frac{x_{k+1} - x_{k}}{\tau} - F^{\top}(y_{k+1} - y_{k}), x_{k+1} - x^{\star} \right\rangle  + \left\langle \frac{y_{k+1} - y_{k}}{\sigma} - F(x_{k+1} - x_{k}), y_{k+1} - y^{\star} \right\rangle \nonumber \\
& \mathrel{\phantom{=}} - \underbrace{\left( \frac{1}{2\tau} \|x_{k+1} - x_{k}\|^2 + \frac{1}{2\sigma} \|y_{k+1} - y_{k}\|^2 - \left\langle F(x_{k+1} - x_{k}), y_{k+1} - y_{k} \right\rangle \right)}_{\mathbf{NE}}, \label{eqn: iter-diff-1-strong-strong}
\end{align}
where $\mathbf{NE}$ represents the numerical error resulting from the implicit discretization. Utilizing the properties of objective functions $f$ and $g^{\star}$, we can further refine the estimation of the iterative difference in~\eqref{eqn: iter-diff-1-strong-strong}. This leads us to the following lemma, which formalizes our findings. 
\begin{lemma}
\label{lem: strong-strong}
Let $f \in \mathcal{S}_{\mu}(\mathbb{R}^{d_1})$ and $g^{\star} \in \mathcal{S}_{\gamma}(\mathbb{R}^{d_2})$. For any step size $\tau\sigma = s^2 \in (0, \|F\|^{-2})$, the discrete Lyapunov function~\eqref{eqn: pdhg-ex-lyapunov} satisfy the following inequality:
\begin{equation}
\label{eqn: lyp-iter-strong-strong}
\mathcal{E}(k+1) - \mathcal{E}(k) \leq - \left( \mu \|x_{k+1} - x^{\star}\|^2 + \gamma \|y_{k+1} - y^{\star}\|^2 \right). 
\end{equation}
\end{lemma}

\begin{proof}[Proof of~\Cref{lem: strong-strong}]
Given that the step sizes $\tau$ and $\sigma$ satisfy the relationship $\tau\sigma = s^2 \in (0, \|F\|^{-2})$, by applying the Cauchy-Schwarz inequality, we can establish the numerical error $\mathbf{NE} \geq 0$. Hence, we can estimate the iterative difference~\eqref{eqn: iter-diff-1-strong-strong} as: 
\begin{align}
& \mathcal{E}(k+1)  - \mathcal{E}(k)  \nonumber  \\
& = \left\langle \frac{x_{k+1} - x_{k}}{\tau} - F^{\top}(y_{k+1} - y_{k}), x_{k+1} - x^{\star} \right\rangle  + \left\langle \frac{y_{k+1} - y_{k}}{\sigma} - F(x_{k+1} - x_{k}), y_{k+1} - y^{\star} \right\rangle. \label{eqn: iter-diff-2-strong-strong}
\end{align}
The optimal conditions of  \eqref{eqn: pdhg-descent} and \eqref{eqn: pdhg-ascent} imply that
\begin{subequations}
\begin{empheq}[left=\empheqlbrace]{align}
  & \frac{x_{k+1} - x_{k}}{\tau}      - F^{\top} (y_{k+1} - y_{k}) + F^{\top} y_{k+1} + \partial f(x_{k+1}) \ni 0,        \label{eqn: pdhga1-descent-ss} \\ 
  & \frac{y_{k+1} - y_{k}}{\sigma} - F (x_{k+1} - x_{k})   -  Fx_{k+1}   + \partial g^{\star}(y_{k+1}) \ni 0.                \label{eqn: pdhga1-ascent-ss}
\end{empheq}
\end{subequations}
Since $f$ and $g$ are strongly convex, we can derive from \eqref{eqn: pdhga1-descent-ss} and \eqref{eqn: pdhga1-ascent-ss} that
\begin{equation*}
	f(x^\star)-f(x^{k+1})\geq \left<\frac{x_{k+1} - x_{k}}{\tau}      - F^{\top} (y_{k+1} - y_{k}) + F^{\top} y_{k+1},x_{k+1}-x^\star\right>+\frac{\mu}{2}\|x_{k+1}-x^\star\|^2,
\end{equation*}
and
$$g^\star(y^\star)-g^\star(y_{k+1})\geq\left\langle \frac{y_{k+1} - y_{k}}{\sigma} - F (x_{k+1} - x_{k})   -  Fx_{k+1} ,y_{k+1}-y^\star\right\rangle+\frac{\gamma}{2}\|y_{k+1}-y^\star\|^2.$$
Similarly, by the properties of saddle point $(x^\star,y^\star)$ in \eqref{eqn: saddle-x-defn}-\eqref{eqn: saddle-y-defn}, we have 
\begin{equation*}
	f(x^\star)-f(x^{k+1})\leq \left<F^\top y^\star,x_{k+1}-x^\star\right>-\frac{\mu}{2}\|x_{k+1}-x^\star\|^2,
\end{equation*}
and
$$g^\star(y^\star)-g^\star(y_{k+1})\leq-\left\langle Fx^\star ,y_{k+1}-y^\star\right\rangle-\frac{\gamma}{2}\|y_{k+1}-y^\star\|^2.$$
Combining the last four inequalities, we deduce that
\begin{subequations}
\begin{empheq}[left=\empheqlbrace]{align}
  & \left\langle \frac{x_{k+1} - x_{k}}{\tau} -  F^{\top}(y_{k+1} - y_{k})   +  F^{\top}(y_{k+1} -y^{\star}), x_{k+1} - x^{\star} \right\rangle \leq -\mu \|x_{k+1} - x^{\star}\|^2,       \label{eqn: pdhga2-inq1-ss} \\ 
  & \left\langle \frac{y_{k+1} - y_{k}}{\sigma} -  F (x_{k+1} - x_{k}) - F(x_{k+1} -x^{\star}), y_{k+1} - y^{\star}  \right\rangle \leq -\gamma \|y_{k+1} - y^{\star}\|^2.                                \label{eqn: pdhga2-inq2-ss}
\end{empheq}
\end{subequations}
By summing the inequalities~\eqref{eqn: pdhga2-inq1-ss} and~\eqref{eqn: pdhga2-inq2-ss} and substituting the result into the iterative difference~\eqref{eqn: iter-diff-2-strong-strong}, we complete the proof. 
\end{proof}

Using the Cauchy-Schwarz inequality, we can estimate the discrete Lyapunov function~\eqref{eqn: pdhg-ex-lyapunov} as follows:
\begin{equation}
\label{eqn: estimate-lyapunov-ss}
\mathcal{E}(k+1) \leq \frac{1+s\|F\|}{2} \cdot \left( \frac{\|x_{k+1} - x^{\star}\|^2}{\tau} + \frac{ \|y_{k+1} - y^{\star}\|^2}{\sigma} \right). 
\end{equation}
This inequality serves as an upper bound for $\mathcal{E}(k+1)$.  Consequently, by comparing it with the right side of  the inequality~\eqref{eqn: lyp-iter-strong-strong} and appropriately scaling the corresponding terms, we can estimate the following proportion:
\begin{equation}
\label{eqn: strong-strong-1}
\frac{\mathcal{E}(k+1) - \mathcal{E}(k)}{\mathcal{E}(k+1)} \leq - \frac{2}{1+s\|F\|} \cdot \min\left\{ \mu \tau, \gamma \sigma \right\}.
\end{equation}
Given that the step sizes satisfy $\tau\sigma = s^2$, we can deduce from inequality~\eqref{eqn: strong-strong-1} the optimal choice of $\tau$ and $\sigma$:
\begin{equation}
\label{eqn: strong-strong-tau-sigma}
\tau = s\sqrt{\frac{\gamma}{\mu}}, \qquad \sigma = s\sqrt{\frac{\mu}{\gamma}},
\end{equation}
which takes the maximum value of $\min\left\{ \mu \tau, \gamma \sigma \right\}$. Taking the step sizes in \eqref{eqn: strong-strong-tau-sigma}, it is derived for the iterative inequality of the discrete Lyapunov function \eqref{eqn: pdhg-ex-lyapunov} as
$$\mathcal{E}(k+1)\leq\frac{1+s\|F\|}{1+s\|F\|+2s\sqrt{\mu\gamma}}\mathcal{E}(k).$$
Hence, we conclude this section with the following theorem.

\begin{theorem}
\label{thm: strong-strong}
By taking the step sizes as given in~\eqref{eqn: strong-strong-tau-sigma} and under the same assumption of~\Cref{lem: strong-strong}, the discrete Lyapunov function~\eqref{eqn: pdhg-ex-lyapunov} holds for any $k \geq 0$ that: 
\begin{equation}
\label{eqn: lyapunov-estimate-1-ss}
\mathcal{E}(k) \leq \left(\frac{1 + s\|F\|}{1+s\|F\| + 2s\sqrt{\mu\gamma}} \right)^{k} \mathcal{E}(0). 
\end{equation}
Furthermore, the iterative sequence $\{(x_k, y_{k})\}_{k=0}^{\infty}$ generated by the~\texttt{PDHG} algorithm,~\eqref{eqn: pdhg-descent} -~\eqref{eqn: pdhg-ascent}, converges to $(x^{\star}, y^{\star})$ at the rate characterized by 
\begin{multline}
\label{eqn: final-estimate-1-ss}
\mu \|x_k - x^{\star}\|^2 + \gamma \|y_k - y^{\star}\|^2 \\ \leq \frac{1+ s\|F\|}{1 - s\|F\|} \left(\frac{1 + s\|F\|}{1+s\|F\| + 2s\sqrt{\mu\gamma}} \right)^{k}  \left(\mu \|x_0 - x^{\star}\|^2 + \gamma \|y_0 - y^{\star}\|^2 \right). 
\end{multline}
\end{theorem}

\section{Conclusion and discussion}
\label{sec: conclu}

In this paper, we harness the power of discrete Lyapunov functions, which can accommodate multiple scales of step sizes, to extend the Lyapunov analysis presented in \citep{li2024understanding} from fixed step sizes to iteration-varying step sizes. By exploiting the nuances of objective functions and saddle points, we appropriately choose a family of step sizes,  $\tau_k \propto (k+1)^{-1}$ and $\sigma_k \propto (k+1)$, ensuring that the~\texttt{PDHG} algorithm converges at a rate approaching $O(1/k^2)$.  Furthermore, in the specific scenario where $\tau_{k+1}\sigma_k = s^2$ and $\theta_k = \tau_{k+1}/\tau_k \in (0, 1)$ as suggested in~\citep{chambolle2011first}, an even faster convergence rate of $O(1/k^2)$ can be achieved. To validate these results, we design a novel discrete Lyapunov function notable for its clarity and simplicity, which elegantly demonstrates the enhanced convergence properties of the \texttt{PDHG} algorithm under these specified conditions. Finally, utilizing this discrete Lyapunov function, we establish the optimal linear convergence rate in the scenario where both the objective functions exhibit strong convexity. 

Recent advancements in these fields have significantly enhanced our understanding of the convergence behaviors of the vanilla gradient descent and illuminated the mechanisms behind accelerated algorithms. Through the application of Lyapunov analysis and phase-space representation, a comprehensive framework grounded on high-resolution ODEs has been well established in the studies~\citep{shi2022understanding, shi2019acceleration, chen2022gradient, chen2022revisiting, chen2023underdamped, li2024linear}.  This framework has been extended to encompass proximal algorithms, like the~\textit{iterative
shrinkage thresholding algorithm} (\texttt{ISTA}) and  the~\textit{fast iterative shrinkage thresholding algorithm}~\texttt{FISTA}s, as detailed in~\citep{li2022linear, li2022proximal, li2024linear}.  Recently, this framework has been further expanded to encompass more advanced algorithms such as the~\textit{Alternating Direction Method of Multipliers}, analyzed through a novel system of high-resolution ODEs in~\citep{li2024understanding1}.  Looking ahead,  there is a compelling opportunity to extend this framework to practical implmentations of~\texttt{ADMM}, exploring methodologies such as linearization~\citep{chambolle2016introduction} and full Jacobian decomposition for multi-separation~\citep{he2015full}.  Additionally, an intriguing avenue for further research lies in investigating traditional minimax optimization algorithms, such as the~\textit{optimal gradient descent ascent} (\texttt{OGDA}) algorithm~\citep{popov1980modification} and the \textit{Extragradient} method~\citep{korpelevich1976extragradient}, and their relationship with \texttt{PDHG} using the insights provided by high-resolution ODEs.

%besides the~\texttt{PDHG} algorithm,  using the , investigate the traditional algorithms, such as \textit{optimal gradient descent ascent} (\texttt{OGDA}) algorithm~\citep{popov1980modification} and the extragradient method~\citep{korpelevich1976extragradient}, their relationship, is also very interesting. 

\section*{Acknowledgements}
This project was partially supported by Grant No.YSBR-034 of CAS, and by the Qingdao Natural Science Foundation, China (No. 23-2-1-158-zyyd-jch).

\bibliographystyle{abbrvnat}
\bibliography{sigproc}
\end{document}